\documentclass[12pt]{article}
\usepackage{amsmath,amssymb}
\newcommand{\G}{\Gamma}
\newcommand{\nin}{\noindent}
\newcommand{\vs}{\vspace*}

\title {\bf The Structures of Zero-divisor Semigroups with Graph $K_n+1$  }
\author{ Tongsuo Wu\thanks{\small Corresponding author. email: tswu@sjtu.edu.cn}\\
Department of Mathematics\\   Shanghai Jiaotong University\\
Shanghai 200030, P. R. China \\
Fan Cheng\thanks {\small email:
chengfan@sjtu.edu.cn}\\
Department of Computer Science and Engineering\\  Shanghai Jiaotong University\\
Shanghai 200030, P. R. China}
\date{}
\begin{document}
\baselineskip=16pt \maketitle \vspace{4mm}
\begin{center}
\begin{minipage}{13cm}

\vs{3mm}\begin{center} {\bf Abstract}\end{center}

\nin{\small   In this paper, we determine the structures of
zero-divisor semigroups whose graph is $K_n+1$, the complete graph
$K_n$ together with an end vertex. We also present a formula to
calculate the number of non-isomorphic zero-divisor semigroups
corresponding to the complete graph $K_n$, for all positive integer
$n$. }

\vs{3mm}\nin {\small\it Key Words:} {\small Commutative zero-divisor
semigroup, complete graph, complete graph with one end vertex }

\vs{3mm}\nin {\small\it 1991 Mathematics Subject Classification:}
{\small 20M14, 05C90}
\end{minipage}
\end{center}

\vs{4mm}\begin{center}{\bf 1. Introduction}\end{center}

\vs{3mm} In this paper, we continue the work in \cite{WULU} of
studying semigroups determined by some graph $G$. In \cite{WULU} it
was proved that the complete graph $K_n$ together with two end
vertices has a unique corresponding zero-divisor semigroup, while
the graph $K_n$ together with more than two end vertices has no
corresponding semigroups, for all $n\ge 4$. In \cite{FL} and
\cite{ZWU}, it was pointed out that both $K_n$ and $K_n$ together
with one end vertex each has multiple corresponding zero-divisor
semigroups. In this paper, we show a formula to calculate the number
of non-isomorphic zero-divisor semigroups corresponding to the
complete graph $K_n$, for all positive integer $n$. We determine the
structure of zero-divisor semigroups whose graph is the complete
graph $K_n$ together with an end vertex. In fact, we have our
discussions according to the four possible values of square of  the
end vertex $x_1$. In three cases, we obtain a simple formula for
counting the number of mutually non-isomorphic zero-divisor
semigroups corresponding to the complete graph $K_n+1$. In the
fourth case (i.e., $x_1^2=x_1$), we give a simple necessary and
sufficient condition and, we give a procedure for listing all the
non-isomorphic zero-divisor semigroups corresponding to the complete
graph $K_n+1$.

For any semigroup $S$, following \cite{B,AL,FR}, associate to $S$ a
simple connected graph $\G(S)$ whose vertex set is $T-\{0\}$, where
$T=Z(S)$ is the set of all zero-divisors of $S$, with $x\not=y$
connected by an edge if $xy=0$. Notice that $T$ is an ideal of $S$
and in particular, it is also a semigroup with the property that it
consists of all zero-divisors of the semigroup $T$. We call such
semigroups $T$ {\it zero-divisor semigroups}. Obviously we have
$\G(S)\cong \G(T)$. For a given connected simple graph $G$, if there
exists a zero-divisor semigroup $S$ such that $\G(S)\cong G$, then
we say that {\it $G$ has corresponding semigroups}, and we call $S$
a {\it semigroup determined by the graph $G$}. Other kinds of
zero-divisor structure were studied in \cite{ZWU, AH}

All semigroups in this paper are multiplicative commutative
zero-divisor semigroups with zero element $0$, where $0x=0$ for all
$x\in S$, and all graphs in this paper are undirected simple and
connected. Throughout this paper, we assume $n\ge 3$.

\vs{4mm}\begin{center} {\bf 2. The graph $K_n$ } \end{center}

\nin\vs{3mm}{\bf Theorem 2.1.} {\it For any n, denote
$M_n=\{0,a_1,a_2,\cdots,a_n\}.$ Then $M_n$ is a zero-divisor
semigroup corresponding to the complete graph $K_n$, if and only if
$M_n$ satisfies the following two conditions:

$(1)$ $a_ia_j = 0, \forall~1\le i\not=j\le n;$

$(2)$ $a_i^{2} = 0$ or $ a_i^{2} = a_i $, or
 $a_i^{2} = a_j$ for some $j \neq i$. In the third case, $a_j^{2} = 0$.}

\vs{3mm}\nin{\bf Proof.} We need only to prove the sufficiency part.
By (1), we only need to check the associative law, namely:
$$(a_ia_j)a_k=a_i(a_ja_k), ~~~~\forall 1\le i ,j ,k \le n. ~~~~ (*)$$
{\bf Case 1}. If $a_i ,a_j,a_k$ are all the same or they are
pairwise different, the associative law obviously holds.

\nin{\bf Case 2}. If $i = j, j \neq k$ ,the right hand of $(*)$ is 0
,and the left of $(*)$ is $a_i^{2}a_k$, if $a_i^{2} \neq a_k$, the
left is 0 too; if $a_i^{2}=a_k$, from (2), we can obtain
$a_k^{2}=0$. Thus in this case the equality $(*)$ holds.

\nin{\bf Case 3}. If $i\neq j, i  = k$, then the left side of $(*)$
is 0, and the right side is 0 too.

\nin{\bf Case 4}. If $i \neq j, j = k$, then the left hand of $(*)$
is 0, while the right side of $(*)$ is $a_ia_j^{2}$. If further
$a_j^{2} \neq a_i$, then the right side is 0. If $a_j^{2} = a_i$,
from (2) we again know $a_i^{2}=0$, so the right hand is 0 too. This
completes the proof of Theorem 2.1.\hfill $\Box$

\vs{3mm}In the following Theorem 2.2, we denote by $p(j,i)$ the
number of the following partitions of the integer $j$: $$d_1  + d_2
+ \cdots + d_i = j,$$where $1 \leq d_1 \leq d_2 \leq \cdots \leq
d_i$.

\vs{3mm}\nin{\bf Theorem 2.2.} {\it The number of non-isomorphic
zero-divisor semigroups corresponding to the complete graph $K_n$ is
$\sum\limits_{k=1}^{n}\sum\limits_{t=0}^{n-k} p(n-t,k) + 1$.}

\vs{3mm}\nin{\bf Proof.} Since $K_n$ is a complete graph, we have
$a_ia_j=0$ for all $i\neq j$, so we only need to decide the value of
$a_i^2$. From Theorem 2.1 we can decompose the set $\{a_i\,|\, 1 \le
i \le n\}$ into a union of the following three pairwise disjoint
subsets:

(1) $A=\{a_i:a_i^2=0 \}$;

(2) $B=\{a_i:a_i^2=a_i \}$;

(3) $C=\{a_i:a_i^2=a_j,a_j \in A \}.$

We assume that the cardinality of $A$ (respectively, $B$) is $|A| =
k\, (0 \leq k \leq n)$ (respectively, $|B| = t\, (0 \leq t \leq
n-k)$). Without loss of generality, we assume the elements in $A$ is
$a_1,a_2,\cdots ,a_k$. From Theorem 2.1 we know that $r^2 \in A$,
$\forall~ r \in C,$ so we can obtain if $k = 0$ , then $t = n$. We
assume there are $\lambda_i$ elements $r\in C$ such that $r^2 =
a_i(1 \leq i \leq k)$, then we get an equation:
    $$\lambda_1 +\lambda_2+ \cdots + \lambda_k = n - t - k.$$
Let $\mu_1,\cdots,\mu_k$ be a permutation of $\lambda_1, \cdots
,\lambda_k$ which satisfies $\mu_1\leq\mu_2\leq\cdots\leq\mu_k$. For
two zero-divisor semigroups $S_1,S_2$ whose zero-divisor graphs are
$K_n$, it is not difficult to see that $S_1$ is isomorphic to $S_2$
$if$ and $only$ $if$ they have the same cardinalities
$|A|,|B|,|C|$,and the same permutation $\mu_1,\cdots,\mu_k$. Thus in
the following we assume that $0 \leq \lambda_1 \leq \cdots \leq
\lambda_k$. So the number of solutions of the equation
$$\lambda_1 +\lambda_2+ \cdots + \lambda_k = n - t - k$$
where $0 \leq \lambda_1 \leq \cdots \leq \lambda_k$, is the number
of corresponding isomorphic zero-divisor semigroups in case $|A|=k
,|B|=t$. Now we substitute $\lambda_i$ by $d_i-1$ $(1 \leq i \leq
k)$, then the above equation is equivalent to :
 $$d_1 + \cdots + d_k  = n - t~~~~(**)$$
where $1 \leq d_1 \leq \cdots \leq d_k , 0 \leq t \leq n-k.$

Finally, we denote by $p(n-t,k)$ the number of solutions of the
equation above, then the number of zero-divisor semigroups
corresponding to the complete graph $K_n$ is
$\sum\limits_{k=1}^{n}\sum\limits_{t=0}^{n-k} p(n-t,k) + 1$.\hfill
$\Box$

\vs{3mm} For any $n\ge 1$, denote
$$s(n)=\sum\limits_{k=1}^{n}\sum\limits_{t=0}^{n-k} p(n-t,k) + 1.$$
One can apply Theorems 2.1 and 2.2 to list all of the twelve (seven)
non-isomorphic zero-divisor semigroups corresponding to $K_4$
($K_3$, respectively). Thus $s(3)=7, s(4)=12.$ These results will be
applied in the last part of the next section.

\vs{4mm}\begin{center} {\bf 3. The graph $K_n$ with one end
vertex}\end{center}

\vs{3mm} Throughout this section, let
$$M_n=\{0,a_1,\cdots, a_n\},~~M_{n,1}=\{0,x_1,a_1,\cdots, a_n\}.$$
and we assume that $M_n$ is a semigroup with $\G(M_n)\cong K_n$. If
$M_{n,1}$ is a commutative zero-divisor semigroup with
$\G(M_{n,1})\cong K_n+1$, {\it the complete graph $K_n$ together
with an end vertex}, then we always assume $a_1x_1=0$ and $x_1$ is
an end vertex. In this case, $M_n$ is an ideal of the semigroup
$M_{n,1}$ by Theorem 4 of \cite{FL}. Thus we have the following
necessary requirements for $M_{n,1}$:

(1) For any $2\le i\le n$, $ a_i^2=0$, or $a_i^2=a_i,$ or for some
$j\not=i$, $a_i^2=a_j$. In the case of $a_i^2=a_j$, we also have
$a_j^2=0$.

(2) $a_1^2=a_1$ or $a_1^2=0$.

(3) $a_1x_1=0$, $a_ix_1\in M_n-\{0\},\forall i\not=1$.

(4) $x_1^2=0,x_1$, or $x_1^2=a_i, i=1,2,\cdots,n$. By symmetry,
one need only consider the four cases of $x_1^2=0,x_1,a_1,a_2$.

In this section, we completely determine the structure of $M_{n,1}$
whose zero-divisor graph is $K_n+1$, the complete graph $K_n$
together with one end vertex. We have our discussions according to
the possible value of $x_1^2$. We remark that for distinct values of
$x_1^2$, the corresponding semigroups are not isomorphic.

\vs{3mm}\nin{\bf Theorem 3.1.} {\it Suppose in $M_{n,1}$ there is a
multiplication such that $a_1x_1=0,a_ia_j=0,\forall i\not=j$. Assume
further $x_1^2=0$. Then $M_{n,1}$ is a semigroup whose zero-divisor
graph is $K_n$ together with an end vertex, if and only if the
following conditions hold:}

(1) $a_1^2=0$.

(2) {\it For all $i\ge 2$, $a_ix_1=a_1$}.

(3) {\it For all $i\ge 2$, $a_i^2=0$ or $a_i^2=a_1$}.

{\it In this situation, there are totally $n$ mutually
non-isomorphic commutative semigroups corresponding to the graph
$K_n+1$}.

\vs{3mm}\nin{\bf Proof.} $\Longrightarrow.$ Suppose $x_1^2=0,
a_1x_1=0$ and assume $M_{n,1}$ is a commutative semigroup such that
$\Gamma(M_{n,1})\cong K_n+1$. For any $i\ge 2$. since
$(a_ix_1)x_1=a_i(x_1^2)=0$, we have $a_ix_1\in \{a_1,x_1\}\cap M_n$.
Thus $a_ix_1=a_1$, $a_i^2x_1=0$. Hence $a_1^2=0$, and $a_i^2\in
\{0,a_1\}$.

$\Longleftarrow.$ Consider the following associative law:
$$(uv)w=u(vw), ~~~~\forall ~u,v,w\in M_{n,1}-\{0\} ~~~~ (*)$$
\nin If $x_1$ does not occur in $u,v,w$, then (*) holds by Theorem
2.1. If $u=v=w=x_1$, then (*) also obviously holds. If exactly two
of $u,v,w$ are the $x_1$, then we have $(x_1v)x_1=x_1(vx_1)$, and
$0=x_1^2w=x_1(x_1w)$ for $w\in M_n$ since either $x_1a_1=0$ or
$x_1a_i=a_1$. In the following we assume that exactly one of the
$u,v,w$ is $x_1$.

Case 1. Assume $u=x_1$ (or equivalently, $w=x_1$). In this case,
$0=(x_1a_1)w$, while $x_1(a_1w)=x_10=0$.
$(x_1a_j)a_i=a_1a_i=x_1(a_ja_i)$ holds for all $i\not=j$. For $i=j$,
we have $(x_1a_i)a_i=0=x_1(a_i^2)$. Thus the multiplication is
associative in this case.

Case 2. Assume $v=x_1$. If $u=a_1$, then
$(a_1x_1)a_k=0=a_1(x_1a_k)$. If $u=a_2$, then
$(a_2x_1)a_k=a_1a_k=0=a_2(x_1a_k)$.

The above discussions show that $M_{n,1}$ is a commutative
zero-divisor semigroup, and $\G(M_{n,1})\cong K_n+1$ if the
conditions (1) to (3) hold. Finally, all the mutually non-isomorphic
commutative zero-divisor semigroups corresponding to $K_n+1$ are
listed in the following
$$M_{n,1}^i=\{0,x_1,a_1,\cdots, a_n\}, ~~i=1,2,\cdots,n,$$
where $a_1x_1=0, x_1^2=0, a_ra_s=0,\forall 1\le r\not=s\le n$,
$a_kx_1=a_1,\forall 2\le k\le n$, and $a_j^2=a_1,\forall 1\le j\le
i$ while $a_j^2=0, \forall i<j\le n$.

This completes the proof. \hfill $\Box$

\vs{3mm}\nin{\bf Theorem 3.2.} {\it Suppose in $M_{n,1}$ there is a
multiplication such that $a_1x_1=0,a_ia_j=0,\forall i\not=j$. Assume
further $x_1^2=x_1$. Then $M_{n,1}$ is a semigroup whose
zero-divisor graph is $K_n+1$, if and only if the following
conditions hold}:

(1) {\it For all $i\ge 2$, $a_ix_1\in M_n-\{a_1, 0\}$ and, there
exists at least one $i\ge 2$ such that $a_ix_1=a_i$}.

(2) {If $a_j=a_ix_1$ ($2\le j\not= i\le n$), then $a_jx_1=a_j$,
$a_j^2=0$, and $a_i^2=0$ or $a_i^2=a_1$};

(3) {\it If $a_rx_1=a_r$ ($r\ge 2$), then $a_r^2$ is equal to one of
the following: $0,a_r,a_j$, where $2\le j\not=r\le n$. If
$a_rx_1=a_r$ and $a_r^2=a_j$ for some $2\le j\not=r\le n$, then
$a_jx_1=a_j, a_j^2=0$}.

(4) {\it $a_1^2=0$ or $a_1^2=a_1$. If $a_i^2=a_1$ for some $i\ge 2$,
then $a_1^2=0$.}

\vs{3mm}\nin{\bf Proof.} $\Longrightarrow.$ Suppose $M_{n,1}$ is a
commutative zero-divisor semigroup such that $\Gamma(M_{n,1})\cong
K_n+1$. Since $n\ge 3$, $x_1$ is an end vertex. Then from
$0=(x_1a_1)a_1=x_1(a_1^2)$ we obtain (4), i.e., either $a_1^2=a_1$
or $a_1^2=0$.

For any $i\ge 2$, we have $a_ix_1\not=a_1$, by the assumption
$x_1^2=x_1$. Thus $a_ix_1\in M_n-\{a_1, 0\}$ since $M_n$ is an ideal
of $M_{n,1}$. This proves the first part of (1). The second
statement of (1) follows from (2). (2) follows easily from the
conditions given.

If $a_rx_1=a_r$ and $a_r^2\not=0,a_r$, then $r\ge 2$ and
$a_r^2=a_j$, where $1\le j\not=r\le n$. If $j=1$, then we have
$0=a_1x_1=a_r^2x_1=a_r^2$, a contradiction. This proves (3).

$\Longleftarrow.$ We only need to check the equality
$$(uv)w=u(vw)~~~(*)$$ for all $u,v,w\in M_{n,1}$.

{\bf Case 1.} If $x_1$ does not occur in $u,v,w$, then (*) holds by
Theorem 2.1. If $u=v=w=x_1$, then (*) also obviously holds.

{\bf Case 2.} Assume that exactly two of $u,v,w$ are the $x_1$. Then
we have $(x_1v)x_1=x_1(vx_1)$. The only other case to verify is
$(x_1x_1)a_i=x_1(x_1a_i)$, i.e., $x_1a_i=x_1(x_1a_i)$: If $i=1$,
then both sides equal to $0$. If $i\ge 2$ and $a_ix_1=a_i$, then
both sides equal to $a_i$. If $i\ge 2$ and $a_ix_1=a_j$ for some
$j\not=i$, then $a_ix_1=a_j=a_jx_1=x_1(x_1a_i)$.

{\bf Case 3.} Now assume that exactly one of the $u,v,w$ is $x_1$.
We need only check in the following two situations.

{\bf Subcase 3.1.} Consider $(x_1v)w=x_1(vw)$. If $v=a_1$, then both
sides equal $0$ since by condition (4), $a_1^2=0~or~a_1$.

If $v=a_i$ ($i\ge 2$) and $a_ix_1=a_i$, then
$(x_1a_i)a_k=a_ia_k=x_1(a_ia_k)$: If $i\not=k$, then each side is
equal to $0$. If $i=k$, then $a_i^2=x_1(a_i^2)$ since $a_i^2$ is
equal to one of the following $0,a_i,a_j$ ($j\not=i$), by condition
(3).

If $v=a_i$ ($i\ge 2$) and $a_ix_1=a_j$ for some $j\not=i$, then the
left side is $(x_1a_i)a_k=a_ja_k$, while the right side is
$x_1(a_ia_k)$. When $j=k$, Then each side is equal to $0$. When
$j\not=k$, again each side is equal to $0$ since $a_i^2=0$ or $a_1$
under assumption $a_ix_1=a_j$ ($i\not=j$).

{\bf Subcase 3.2.} Finally, let us consider
$$(a_ix_1)a_k=a_i(x_1a_k)~~~~~(**)$$ It is easy to verify
$(a_1x_1)a_k=a_1(x_1a_k)$ for all $k$. In the following we assume
$i\ge 2$. If $a_ix_1=a_i$, then the left side of $(**)$ is $a_ia_k$
and the right side is $a_i(x_1a_k)$. If further $i=k$, then both
sides are $a_i^2$. If $i\not=k$, then both sides are $0$. Finally,
we assume $a_ix_1=a_j$ for some $j\not=i$. Then the left side of
$(**)$ is $a_ja_k$ and the right side is $a_i(x_1a_k)$. If
$j\not=k$, then each side is equal to $0$ since $i\ge 2$, and
$a_ix_i\not=x_i$. If $j=k$, then the left side is $a_j^2=0$ and the
right side is $a_i(x_1a_j)=a_ia_j=0$. This completes the whole
verification. \hfill $\Box$

\vs{3mm}\nin{\bf Theorem 3.3.} {\it Suppose in $M_{n,1}$ there is a
multiplication such that $a_1x_1=0,a_ia_j=0,\forall i\not=j$.}

(i) {\it If in addition $x_1^2=a_1$, then $M_{n,1}$ is a semigroup
whose zero-divisor graph is $K_n$ together with an end vertex, if
and only if $a_rx_1=a_1, \forall r\ge 2$ and $a_i^2=0$ for all $i$.
In this situation, there is exactly one zero-divisor semigroup $S$
with graph $\G(S)\cong K_n+1$.}

(ii) {\it If in addition $x_1^2=a_2$, then $M_{n,1}$ is a
semigroup whose zero-divisor graph is $K_n$ together with an end
vertex, if and only if the following conditions hold:}

(1) {\it $a_1^2=0$ and for any $i\ge 3$, $a_ix_1=a_1$ and, $
a_i^2=0$ or $a_i^2=a_1.$}

(2) {\it Exactly one of the following cases occurs}:

\hspace{1cm}(A) {$a_2x_1=a_1$, and $a_2^2=0$};

\hspace{1cm}(B) {\it $a_2x_1=a_2$, and $a_2^2=a_2$};

\hspace{1cm}(C) {\it $a_2x_1=a_r$ for some $3\le r\le n$, and
$a_2^2=a_1,a_r^2=0$.

In the case of (ii), there are totally $3n-4$ mutually
non-isomorphic commutative semigroups corresponding to the graph
$K_n+1$. }

\vs{3mm}\nin{\bf Proof.} (i) Assume that $M_{n,1}$ is a semigroup
whose zero-divisor graph is $K_n$ together with an end vertex. If in
addition $x_1^2=a_1$, then $a_1^2=0$ and for all $i\ge 2$,
$(a_ix_1)x_1=0$. Thus $a_ix_1=a_1, a_i^2=0, \forall i\ge 2$.
Conversely, it is routine to verify that the associative law holds.

(ii) $\Longrightarrow.$  Suppose that $M_{n,1}$ is a semigroup whose
zero-divisor graph is $K_n$ together with an end vertex and assume
$x_1^2=a_2$ Then (1) follows easily from the assumption. Since $M_n$
is an ideal of $M_{n,1}$, thus $a_2x_1\in M_n-\{0\}$. If
$a_2x_1=a_1$, then $a_2^2=(a_2x_1)x_1=a_1x_1=0$. This proves (A). In
a similar manner, one obtains (B) and (C).

$\Longleftarrow.$ Again we need only check the equality
$$(uv)w=u(vw)~~~(*)$$ for all $u,v,w\in M_{n,1}$.

{\bf Case 1.} If $x_1$ does not occur in $u,v,w$, then (*) holds by
Theorem 2.1. If $u=v=w=x_1$, then (*) also obviously holds.

{\bf Case 2.} Assume that exactly two of $u,v,w$ are the $x_1$. Then
we need only to verify $a_2a_i=x_1(x_1a_i)$ since
$(x_1x_1)a_i=a_2a_i$. In fact, if $i=1$, then both are $0$. If $i\ge
3$, then both sides are $0$. If $i=2$, then we need to verify
$a_i^2=x_1(x_1a_2)$. This is the case by the assumption of (A),or
(B), or (C).

{\bf Case 3.}  Assume that exactly one of the $u,v,w$ is $x_1$. Then
we need to verify both $(x_1a_i)a_j=x_1(a_ia_j)$ and
$(a_ix_1)a_j=a_i(x_1a_j)$. In the following we only verify the first
equality because the verifications of the second one is similar.

Consider the possible equality $(x_1a_i)a_j=x_1(a_ia_j)$:

(1) If $i=1$, then each side is equal to $0$ since $a_1^2=0$.

(2) If $i\ge 3$, then the left side is $0$. If in addition, $i=j$,
then the right side is $x_1a_i^2=0$ since $a_i^2=0$ or $a_i^2=a_1$.
If $i\not=j$, then $a_ia_j=0$.

(3) The last subcase is $i=2$. If $a_2x_1=a_1$, then each side is
equal to $0$. If $a_2x_1=a_2$, then
$(x_1a_2)a_2=a_2a_2=a_2=x_1(a_2a_2)$, and
$(x_1a_2)a_j=a_2a_j=0=x_1(a_2a_j)$ for $j\not=2$. If $a_2x_1=a_r$
with $r\ge 3$, then $(x_1a_2)a_2=a_ra_2=0= x_1a_1=x_1(a_2a_2)$, and
for $j\not=2$, $(x_1a_2)a_j=a_ra_j=0= x_10=x_1(a_2a_j)$.

Finally, cases (A) and (B) each has $n-1$ mutually non-isomorphic
commutative semigroups corresponding to the graph $K_n+1$. In case
(C) we have mutually $n-2$ non-isomorphic corresponding commutative
semigroups. This completes the proof. \hfill $\Box$

\vs{3mm}We end up this paper with the following remarks.

\vs{3mm}\nin{\bf Remark 1.} Denote by $k_2(n)$ the number of
mutually non-isomorphic commutative semigroups corresponding to the
graph $K_n+1$ in Theorem 3.2. Then by Theorems 3.1 to 3.3,
$k_2(n)+4n-3$ {\it is the total number of mutually non-isomorphic
commutative semigroups corresponding to the graph $K_n+1$}.

\vs{3mm}\nin {\bf Remark 2.} In the following we provide a procedure
for calculating $k_2(n)$. In Theorem 3.2, let $k_2(n,r)$ be the
number of mutually non-isomorphic commutative semigroups
corresponding to the graph $K_n+1$, in which there are exactly $r$
numbers $i\in \{2,\cdots, n\}$ such that $a_ix_1=a_i$, where
$r=1,2,\cdots,n-1$. Then $$k_2(n)=\sum_{r=1}^{n-1}k_2(n,r).$$

When $r=n-1$, we have $a_ix_1=a_i$ for all $2\le i\le n$. In this
subcase, the value of $a_1^2$ is either $0$ or $a_1$, the value of
$a_i^2$ is $0$, $a_i$ or$a_j$ ($2\le i\not=j\le n$). By Theorem 3.2
and Theorem 2.1, $k_2(n,n-1)=2s(n-1)$, where $s(n-1)$ is the number
of non-isomorphic zero-divisor semigroups corresponding to the
complete graph $K_{n-1}$.

When $r=1$, we can assume $a_2x_1=a_2$. For all $i\ge 3$, we obtain
$a_ix_1=a_2$ by condition (3) of Theorem 3.2 and hence, $a_i^2=0$ or
$a_i^2=a_1$. Furthermore, $a_1^2=a_1$ only if $a_i^2=0$ holds for
all $i\ge 3$. Finally, it is routine to check that there are $n$
mutually non-isomorphic associative multiplication tables in
$M_{n,1}$ such that $\G(M_{n,1})\cong K_n+1$. Hence $k_2(n,1)=n$.

When $r=2$, without loss of generality we can assume $a_ix_1=a_i,
i=2,3$. Then by condition (3) of Theorem 3.2, we have $a_jx_1=a_2,$
or $a_jx_1=a_3$ for all $j\ge 4$. If $\{a_2,a_3\}=\{a_jx_1\,|\,4\le
j\le n\}$, then $a_2^2=0, a_3^2=0$. For any $4\le i\le n$, either
$a_i^2=0$ or $a_i^2=a_1$ by condition (2), and $a_1^2=a_1$ only if
$a_i^2=0,\forall i\ge 4$. Thus in this case, we have $n-1$
multiplication tables. The only other case is
$\{a_2\}=\{a_jx_1\,|\,4\le j\le n\}$. In this case if $n\ge 4$, we
have $3(n-1)$ multiplication tables since the value of $a_3^2$ could
be one of $0,a_2$ or $a_3$. When $r=2, n\ge 5$ and $x_1^2=x_1$,
there are totally $4(n-1)$ mutually non-isomorphic associative
multiplication tables in $M_{n,1}$ such that $\G(M_{n,1})\cong
K_n+1$. Hence
$$k_2(n,2)=\left\{\begin{array}{ll}
3, &\textrm{ if $n=3$}\\
3\times (4-1),& \textrm{ if $n=4$}\\
4(n-1), &\textrm{ if $n\ge 5$.}
\end{array}\right.$$

For $r=3,4,\cdots,n-2$, one can continue these discussions. When
$r=3$, like the $r=2$ case, there are five results for
$n=3,4,5,6,\ge 7$ respectively.

\vs{3mm}If $n$ is small, it is not very difficult to calculate the
number $k_2(n)$. For example, $k_2(3)=6$ and $k_2(4)=4+3\times
(4-1)+2\times 7=27$, $k_2(5)=5+4\times(5-1)+2\times 7+2\times
12=59$. Thus $K_3+1$, $K_4+1$ and $K_5+1$ has $15, 40$ and $76$
mutually non-isomorphic commutative semigroups, respectively. For
general $n$, we still do not know if there is a simple formula for
calculating $k_2(n)$.

\end{document}